\documentclass[12pt]{amsart}

\usepackage{latexsym}
\usepackage{amsmath}
\usepackage{amssymb}
\usepackage{amscd}
\usepackage{amsthm}
\usepackage{xypic}%HC
\xyoption{curve}
\usepackage{ifthen} %HC
\usepackage{hyperref} %HC
\usepackage{graphicx}

\theoremstyle{plain}
 \newtheorem{theorem}{Theorem}[section]
 \newtheorem{proposition}[theorem]{Proposition}
 \newtheorem{lemma}[theorem]{Lemma}

\theoremstyle{definition}
 \newtheorem{definition}[theorem]{Definition}
 \newtheorem*{example}{Example}

\theoremstyle{remark}
 \newtheorem{remark}[theorem]{Remark}

\makeindex

\begin{document}
\title[Rationality of moduli spaces of plane curves]{Rationality of the moduli spaces of plane curves of sufficiently large 
degree} 

\author{Christian B\"ohning\\ 
 Hans-Christian Graf von Bothmer}

\maketitle

\newcommand{\PP}{\mathbb{P}} %HC
\newcommand{\QQ}{\mathbb{Q}} %HC
\newcommand{\ZZ}{\mathbb{Z}} %HC
\newcommand{\CC}{\mathbb{C}} %HC
\newcommand{\rmprec}{\wp}%HC
\newcommand{\rmconst}{\mathrm{const}}%HC
\newcommand{\xycenter}[1]{\begin{center}\mbox{\xymatrix{#1}}\end{center}} %HC

\newboolean{xlabels} %HC
\newcommand{\xlabel}[1]{ %HC
                        \label{#1} %HC
                        \ifthenelse{\boolean{xlabels}} %HC
                                   {\marginpar[\hfill{\tiny #1}]{{\tiny #1}}} %HC
                                   {} %HC
                       } %HC
%\setboolean{xlabels}{true} %HC
\setboolean{xlabels}{false} %HC

\begin{abstract}
We prove that the moduli space
%$C (d) := \PP\, H^0 (\PP^2, \mathcal{O} (d)) / \mathrm{PGL}_3\,\CC$ 
of plane curves of degree $d$ is rational for all sufficiently large $d$.
\end{abstract}

\section{Introduction}

It is a classical question, which can be traced back to works of Hilbert and Emmy Noether, whether the orbit spaces $\mathbb{P}/
G$ are rational where $\mathbb{P}$ is a projective space and $G$ is a reductive algebraic group acting linearly in $\mathbb{P}$.
If $G$ is not assumed connected, in fact for $G$ a finite solvable group, D. Saltman has shown in \cite{Sa} that the answer to
this question is negative in general (Emmy Noether had apparently conjectured that the quotient should be rational in this case).
No counterexamples are known for connected complex reductive groups
$G$.

\

For simply connected classical groups, the quotients $\mathbb{P}/G$ are known to be stably rational by results of F.
Bogomolov
\cite{Bog1}. The question whether stably rational varieties are always rational is the well-known Zariski problem which
Beauville, Colliot-Th\'{e}l\`{e}ne, Sansuc and Swinnerton-Dyer \cite{Beau} answered in the negative as well: There are
three-dimensional conic bundles $X$ over rational surfaces which are irrational, but $X\times\mathbb{P}^3$ is rational. This uses
the method of intermediate Jacobians by Clemens-Griffiths \cite{C-G} which, however, seems to work only for threefolds. In
general, it is rather hard to distinguish stably rational and rational varieties. The method connected with the birational
invariance of the Brauer-Grothendieck group used previously by Artin and Mumford \cite{A-M} to obtain more elementary examples of
unirational non-rational threefolds, is insensitive to this distinction (unirational varieties constitute a strictly
bigger class than stably rational ones; e.g. Saltman's counterexamples mentioned above are not even stably rational). The reader
may find these and other methods to prove irrationality, including the use of Noether-Fano inequalities via untwisting of
birational maps and Koll\'{a}r's method of differential forms in characteristic $p$ to prove non-rationality of some general
hypersurfaces, in the survey by V. A. Iskovskikh and Yu. G. Prokhorov
\cite{Is-Pr}.

\

The geometrically most relevant case of the general question discussed above seems to be the case of the moduli space of
projective hypersurfaces of degree $d$ in $\mathbb{P}^n$, which we denote by $\mathrm{Hyp} (d,\: n)$. Here rationality is known in
the following cases:
\begin{itemize}
\item
$n=1$ (the classical case of binary forms resp. sets of points on the projective line), $d$ odd \cite{Kat}, $d$ even \cite{Bo-Ka},
\cite{Bog2}
\item
$n=2$, $d\le 3$ (well known), $d=4$ (\cite{Kat1}, \cite{Kat2}), $d \equiv 1$ (mod $4$) (\cite{Shep}), $d \equiv 1$ (mod $9$) and
$d\ge 19$ (\cite{Shep}); $d\equiv 0$ (mod $3$) and $d\ge 1821$ \cite{Kat0} (this article contains the remark that the author
obtained the result also for $d\ge 210$, unpublished); the same paper also gives some results for congruences to the modulus $39$;
furthermore, there are some unpublished additional cases in the case of plane curves which we do not try to enumerate.
\item
$n=3$, $d\le 2$ (obvious), $d=3$ (Clebsch and Salmon; but see \cite{Be}).
\item
$n > 3$, $d\le 2$ (obvious).
\end{itemize}
This represents what we could extract from the literature. It is hard to say if it is exhaustive. The reader may consult the very
good (though not recent) survey article \cite{Dol0} for much more information on the rationality problem for fields of
invariants.\\
The main theorem of the present article is

\begin{theorem} \xlabel{tMain}
The moduli space of plane curves of sufficiently large degree $d >>0$ under projective equivalence is rational.
\end{theorem}

More precisely, for $d=3n$, $d\ge 1821$, this was proven by Katsylo \cite{Kat0} as a glance back at the preceding summary shows.
We use this result and don't improve the bound for $d$. For $d \equiv 1$ (mod $3$), we obtain rationality for $d\ge 37$. For $d
\equiv 2$ (mod$3$), we need $d\ge 65$.

\

Let us turn to some open problems. First of all, the method used in this paper seems to generalize and -provided the required
genericity properties hold and can be verified computationally- could yield a proof of the rationality of $\mathrm{Hyp} (d,\: n)$
for
\emph{fixed}
$n$ if the degree $d$ is large enough and $n+1$ does not divide $d$. The latter case might be amenable to the techniques of
\cite{Kat0} in general. Thus the case of the moduli spaces of surfaces of degree $d$ in $\mathbb{P}^3$ seems now tractable with
some diligence and effort. But we do not see how one could obtain results on $\mathrm{Hyp} (d, \: n)$ \emph{for all} $n$ (and
$d$ sufficiently large compared to $n$). 

More importantly, whereas we think that it is highly plausible that $\mathrm{Hyp}(d,\: n)$ is always
rational if
$d$ is sufficiently large compared to $n$, we do not want to hazard any guess in the case where $d$ is small. In fact, we do not
know any truely convincing philosophical reason why $\mathrm{Hyp}(d,\: n)$ should be rational in general; the present
techniques of proving rationality always seem to force one into assuming that $d$ is sufficiently large if one wants to obtain
an infinite series of rational examples by a uniform method. Moreover, it can be quite painstaking and tricky to get a hold of
the situation if $d$ is small as Katsylo's tour de force proof for $\mathfrak{M}_3$ (i.e. $\mathrm{Hyp}(4, \: 2)$) in 
\cite{Kat1}, \cite{Kat2} amply illustrates. The maybe easiest unsolved cases are $\mathrm{Hyp}(6, \: 2)$ (plane sextics) and
$\mathrm{Hyp}(4,\: 3)$ (quartic surfaces). Note that the former space is birational to the moduli space of polarized K3 surfaces
$(S, \: h)$ of degree $2$ (thus $S$ is a nonsingular projective K3 surface and $h\in\mathrm{Pic}(S)$ is the class of an ample
divisor with $h^2=2$), and the latter space is birational to the moduli space of polarized K3 surfaces of degree $4$. 

We would like to thank Fedor Bogomolov and Yuri Tschinkel for suggesting this problem and helpful
discussions. 

% begin HC
\section{outline of Proof}

The structural pattern of the proof is similar to \cite{Shep}; there the so called \emph{method of covariants} is
introduced, and we learnt a lot from studying that source.\\
We fiber the space $\mathbb{P} \left( \mathrm{Sym}^d\, ( \mathbb{C}^3)^{\vee } \right)$ of degree $d$ plane curves over the space
of plane quartics, if $d\equiv 1$ (mod $3$), and over the space of plane octics, if $d\equiv 2$ (mod $3$), i.e. we construct
$\mathrm{SL}_3 (\mathbb{C})$-equivariant maps 
\[
	S_d \colon  \mathrm{Sym}^d\, ( \mathbb{C}^3)^{\vee } \to \mathrm{Sym}^4\, ( \mathbb{C}^3)^{\vee }
\]
and
\[
	T_d \colon  \mathrm{Sym}^d\, ( \mathbb{C}^3)^{\vee } \to \mathrm{Sym}^8\, ( \mathbb{C}^3)^{\vee }
\]
($S_d$ coincides with the covariant used in \cite{Shep} for the case $d=9n+1$). These maps are of degree $4$ as polynomials in the
coordinates on the source spaces, i.e. of degree $4$ in the curve coefficients. They are constructed via the symbolic method recalled
in section 3. Furthermore they induce \emph{dominant} rational maps on the associated projective spaces. We remark
here that the properties of
$S_d$ and
$T_d$ essential for the proof are that they are of fixed low degree in the curve coefficients, take values in spaces of curves of
fixed low degree, and are sufficiently generic.

We now focus on the case $d \equiv 1$ (mod $3$). The proof has three main steps:
\begin{itemize}
\item[(1)]
$\mathrm{Hyp} (4, \: 2)$ is stably rational, more precisely its product with $\mathbb{P}^8$
is rational; cf. \cite{Bo-Ka}, Theorem 1.1 for this.  
\item[(2)]
We find a linear subspace $L_S \subset \mathrm{Sym}^d (\mathbb{C}^3)^{\vee
}$ such that $\PP(L_S)$ is contained in the base locus $B_{S_d}$ of $S_d$ with
a full triple structure, i.e.
$I_{\PP(L_S)}^3
\supset I_{B_{S_d}}$, and consider the projection $\pi_{L_S}$ away from  $\mathbb{P}(L_S)$  onto $\PP (\mathrm{Sym}^d
(\mathbb{C}^3 )^{\vee }/L_S)$. We show that 
 a general fibre of $S_d$  is birationally  a vector bundle over a rational base.
\item[(3)]
The quotient map \[ \mathbb{P} \mathrm{Sym}^4\, ( \mathbb{C}^3)^{\vee } \dasharrow \left( \mathbb{P} \mathrm{Sym}^4\, (
\mathbb{C}^3)^{\vee }\right) /\mathrm{PGL}_3 (\mathbb{C}) \] 
has a section $\sigma_4$. Pulling back the linear fibrations constructed in $(2)$ via $\sigma_4$ we show that the moduli space of plane
curves of degree $d$ is birational to $\mathrm{Hyp} (4, \: 2) \times \mathbb{P}^{N}$, where $N$ is  large, whence we conclude by (1).
\end{itemize}

The main computational difficulty occurs in (2) where we have to establish that $L_S$ is sufficiently generic. Projecting from
$\PP(L_S)$ we obtain a diagram

\xycenter{
	\PP \mathrm{Sym}^d (\mathbb{C}^3 )^{\vee } \ar@{-->}[r]^{S_d} \ar@{-->}[d]^{\pi_{L_S}}& \PP \mathrm{Sym}^4 (\mathbb{C}^3 )^{\vee
}\\
	\PP (\mathrm{Sym}^d (\mathbb{C}^3 )^{\vee }/L_S)
	}

We show that

\begin{itemize}
\item[($\ast $)]
for a particular (hence a general) $\bar{g} \in \PP (\mathrm{Sym}^d (\mathbb{C}^3 )^{\vee }/L_S)$ the map 
\[
	S_d|_{\mathbb{P} (L_S+\mathbb{C} g )} \colon \mathbb{P} (L_S+\mathbb{C} g)  \dasharrow  \PP \mathrm{Sym}^4 (\mathbb{C}^3 )^{\vee
}
\]
is surjective. 
\end{itemize}
Note that $S_d|_{\mathbb{P} (L_S+\mathbb{C} g )}$ is linear since $L_S$ is contained in the
base locus with full triple structure and $S_d$ is of degree $4$ in the curve coefficients. It is therefore
enough to explicitly construct points in the image that span $\PP \mathrm{Sym}^4 (\mathbb{C}^3
)^{\vee }$. From this it follows at once that a general fibre of $S_d$ is mapped dominantly by $\pi_{L_S}$ whence we may view
such a fibre birationally as a vector bundle over a rational base. To understand better why the dominance of $S_d$ is not
sufficient here, it is instructive to keep the following example in mind:

\begin{figure} 
\includegraphics*[width=6cm]{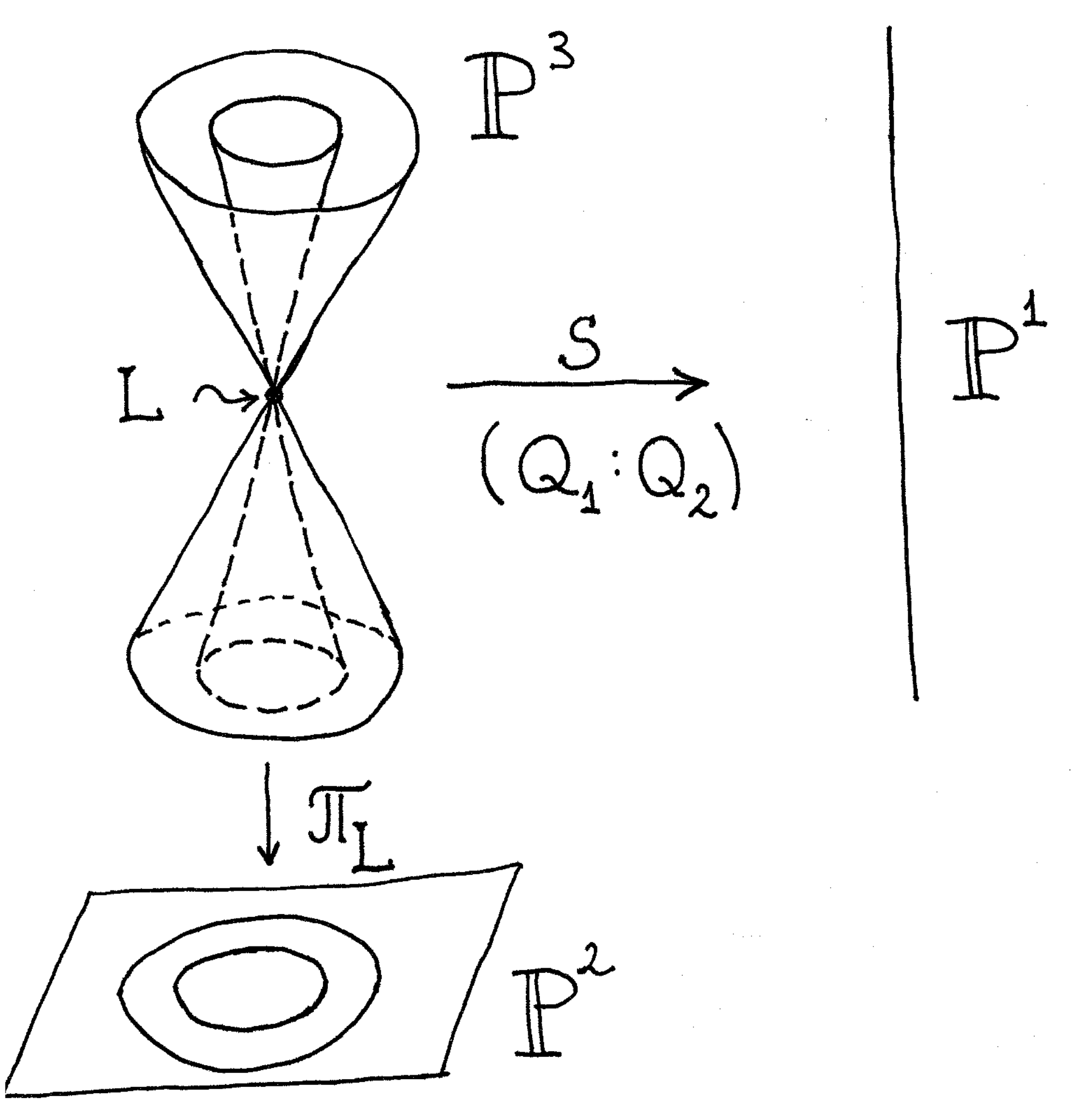}
\caption{A non generic $L$.} \label{fExample}
\end{figure}

\begin{example}
Consider the rational map 
$$
S \colon \PP^3 \dasharrow \PP^1,\quad
	x \mapsto (Q_1(x) : Q_2(x))
$$	
where $Q_1, Q_2$ are quadric cones with vertex $L$. $S$ is dominant. Projection from the vertex
$L$ to $\PP^2$ is also dominant, but the quadric cones (i.e. the fibers of $S_d$) do not map dominantly to $\PP^2$ (see Figure \ref{fExample}). The projection fibers are the lines through $L$ and indeed each such line is contained completely in one cone in the pencil $\lambda Q_1 + \mu Q_2$. 

The base locus $B$ of $S$ consists of $4$ lines that meet in $L$. If on the other hand we project from a smooth point $L' \subset B$ then a general fiber of $S$ maps dominantly to $\PP^2$. Indeed a general line through $L'$ intersects all cones. 
\end{example}

The complications in proving ($\ast )$ arise due to the fact that the natural description of $L_S$ is in terms of the
\emph{monomials} which span it, whereas $S_d$ can be most easily evaluated on forms which are written \emph{as sums of powers of
linear forms}. These two points of view do not match, and we cannot repose on methods in \cite{Shep}.  Instead we introduce new techniques in section \ref{sFiberwise} to solve this difficulty:
\begin{itemize}
\item
We use \emph{interpolation polynomials} to write down elements in $L_S$ as \emph{sums of powers of linear forms}.
\item
Next we employ considerations of \emph{leading terms} (or, geometrically, jets at infinity) to eliminate the interpolation polynomials,  from our formulae. 
\item
For large enough $d=3n+1$, we finally reduce (4) to the property that a certain matrix $M(n)$ has full rank. The size
of $M(n)$ is independent of $n$ while (\emph{and this is the main point}) its entries are of the form
\[
		\sum_\nu \rho_\nu^n P_\nu(n)
\]
where $P_\nu$ are polynomials of fixed degree (i.e. independent of $n$), $\rho_\nu$ are constants, and the number of summands in the
expression is independent of $n$. This is possible only because we eliminated the interpolation polynomials in the previous step. 
\item
By choosing a point
$g$ with integer coefficients we
 can arrange that $\rho_\nu$ and $P_\nu(n)$ are defined over $\QQ$ with denominators that are not divisible by a small prime
$\rmprec$ which we call the precision of our calculation. Thus if we work over \emph{the finite field}
$\mathbb{F}_{\rmprec}$, the matrix
$M(n)$ is periodic in
$n$ with period
$\rmprec(\rmprec-1)$. A computer algebra program is then used to check that these matrices all have full rank. By semicontinuity,
this proves that $M(n)$ has full rank for all $n$ in characteristic $0$.
\end{itemize}
In a rather round-about sense, we have also been guided by the principle that \emph{evaluation of a polynomial} at a special point
can be much cheaper than \emph{computing the polynomial}.

% end HC

\section{Notation and definition of the Covariants}
For definiteness, the base field will be $\CC$, the field of complex numbers, though one might replace it by any
algebraically closed field of characteristic $0$ throughout.\\
Let $G:=\mathrm{SL}_3\, ( \CC )$, and let $\bar{G}:=\mathrm{PGL}_3\, ( \CC )$ be the
adjoint form of
$G$. We denote by
$V(k)$ the irreducible $G$-representation $\mathrm{Sym}^k\, (\CC^3)^{\vee }$. We fix a positive integer $d$ not divisible by $3$, $d=3n+1$ or $d=3n+2$, $n\in
\mathbb{N}$.\\
The symbol $[k]$, $k\in\mathbb{N}$, denotes the set of integers from $0$ (incl.) to $k$ (incl.). Let $x_1,\: x_2, \: x_3\in
(\CC^3)^{\vee }$ denote the basis dual to the standard basis in $\CC^3$ and put $\mathbf{x} := (x_1,\: x_2, \:
x_3)$. We will use Schwartz's multi-index notation and denote multi-indices by lower case boldface letters. Thus we write a
general homogeneous form $f\in V(d)$ of degree $d$ as
\begin{gather}
f = \sum\limits_{\stackrel{\mathbf{i} \in [d]^3}{| \mathbf{i} | = d}} \frac{d !}{\mathbf{i} !} A_{\mathbf{i}}
\mathbf{x}^{\mathbf{i}}\, ,
\end{gather}
where $\mathbf{i} ! = i_1! i_2 ! i_3 !$, $|\mathbf{i} | := i_1+i_2+i_3$, $A_{\mathbf{i}}= A_{i_1i_2i_3}$,
$\mathbf{x}^{\mathbf{i}} = x_1^{i_1}x_2^{i_2}x_3^{i_3}$. We will use the symbolical method introduced by Aronhold and Clebsch
to write down
$G$-equivariant maps (covariants) from $V(d)$ to $V(4)$ (if $d=3n+1$) or to $V(8)$ (if $d=3n+2$). It is explained in \cite{G-Y} and, from a modern point of view, in
\cite{Dol1}, chapter 1. We denote by $\mathbf{\alpha} = (\alpha_1, \: \alpha_2, \:\alpha_3)$ a vector of symbolic variables,
and also introduce vectors $\mathbf{\beta}$, $\mathbf{\gamma}$, $\mathbf{\delta}$, similarly. We write
$\mathbf{\alpha}_{\mathbf{x}} =\alpha_1 x_1 +\alpha_2x_2 +\alpha_3x_3$, and similarly $\mathbf{\beta}_{\mathbf{x}}$,
$\mathbf{\gamma}_{\mathbf{x}}$,
$\mathbf{\delta}_{\mathbf{x}}$. Moreover we define the bracket factor $(\mathbf{\alpha}\,\mathbf{\beta}\, \mathbf{\gamma})$ by
\begin{gather*}
(\mathbf{\alpha}\,\mathbf{\beta}\, \mathbf{\gamma}) := \det \left( \begin{array}{ccc} \alpha_1 & \alpha_2 & \alpha_3\\ \beta_1 &
\beta_2 & \beta_3 \\ \gamma_1 & \gamma_2 & \gamma_3
\end{array}\right)
\end{gather*}
and write $(\mathbf{\alpha}\,\mathbf{\beta}\, \mathbf{\delta})$ etc. similarly. The idea in this calculus is to write
$f\in V(d)$ symbolically as a power of a linear form in several ways:
\begin{gather}
f = \mathbf{\alpha}_{\mathbf{x}}^d = \mathbf{\beta}_{\mathbf{x}}^d =\mathbf{\gamma}_{\mathbf{x}}^d
=\mathbf{\delta}_{\mathbf{x}}^d \, ,
\end{gather}
whence the identities
\begin{gather}
A_{\mathbf{i}} = \mathbf{\alpha}^{\mathbf{i}} = \mathbf{\beta}^{\mathbf{i}} =\mathbf{\gamma}^{\mathbf{i}} =
\mathbf{\delta }^{\mathbf{i}} \, .
\end{gather}
If $d=3n+1$, define a covariant $S_d \, :\, V(d) \to V(4)$ of order $4$ and degree $4$ by the following prescription:
\begin{align}
I(\mathbf{\alpha}, \: \mathbf{\beta},\:\mathbf{\gamma}, \: \mathbf{\delta}) &:= (\mathbf{\alpha}\,\mathbf{\beta}\,
\mathbf{\gamma}) (\mathbf{\alpha}\,\mathbf{\beta}\, \mathbf{\delta}) (\mathbf{\alpha}\,\mathbf{\gamma}\, \mathbf{\delta})
(\mathbf{\beta}\,\mathbf{\gamma}\, \mathbf{\delta})\, \\
S_d (\mathbf{\alpha}, \: \mathbf{\beta},\:\mathbf{\gamma}, \: \mathbf{\delta}) &:= I^n \mathbf{\alpha}_{\mathbf{x}}
\mathbf{\beta}_{\mathbf{x}}  \mathbf{\gamma}_{\mathbf{x}}  \mathbf{\delta}_{\mathbf{x}} \, .
\end{align}
The formula for $S_d$ should be read in the following way: The right hand side of $(5)$, when we multiply it out formally, is a
sum of monomials $\mathbf{\alpha}^{\mathbf{i}} \mathbf{\beta}^{\mathbf{j}} \mathbf{\gamma}^{\mathbf{k}}
\mathbf{\delta}^{\mathbf{l}} \mathbf{x}^{\mathbf{e}}$, $\mathbf{i}, \: \mathbf{j}, \:\mathbf{k}, \: \mathbf{l} \in [d]^3$,
$\mathbf{e}\in [4]^3$, and $|\mathbf{i}| = \dots = |\mathbf{l}| = d$, $|\mathbf{e}| = 4$. Thus one can use equations (1) and
(3) to rewrite the right-hand-side unambiguously in terms of the coefficients $A_{\mathbf{i}}$ of $f\in V(d)$. Hence $S_d$ may
be viewed as a map from $V(d)$ to $V(4)$, homogeneous of degree $4$ in the coefficients $A_{\mathbf{i}}$, which is clearly
$G$-equivariant. By abuse of notation, we denote the induced rational map
\begin{gather}
S_d\, :\, \PP\, V(d) \dasharrow  \PP\, V(4)
\end{gather}
by the same letter. Note that $I$ defined by equation $(4)$ may be viewed as an invariant of plane cubics $I\, : V(3) \to
\CC$ of degree $4$ in the coefficients of the cubic. In fact, this is the famous Clebsch invariant, vanishing on the
locus of Fermat cubics, or vanishing on the equi-anharmonic cubics, i.e. nonsingular plane cubics which can be written as a
double cover of $\PP^1$ branched in four points with equi-anharmonic cross-ratio. Equi-anharmonic cross-ratio means
cross-ratio equal to minus a cube root of $1$. Equi-anharmonic quadruples of points in $\PP^1$ are one of the two
possible
$\mathrm{PGL}_2\,\CC$-orbits of $4$ points in $\PP^1$ with non-trivial isotropy group (the other orbit being
quadruples with harmonic cross-ratio, i.e. equal to $-1$, $1/2$ or $2$). See \cite{D-K}, (5.13), for details.\\
The letter $S$ in $S_d$ was chosen in honor of the 19th century Italian geometer Gaetano Scorza, who studied in detail the
map $S_4$, called the \emph{Scorza map} (cf. \cite{D-K}, \S 6 and \S 7, and \cite{Dol2}, section 6.4.1). 

\

Similarly, for $d=3n+2$, we define a covariant $T_d\, :\, V(d) \to V(8)$ of order $8$ and degree $4$ by
\begin{align}
T_d (\mathbf{\alpha}, \: \mathbf{\beta},\:\mathbf{\gamma}, \: \mathbf{\delta}) &:= I^n \mathbf{\alpha}_{\mathbf{x}}^2
\mathbf{\beta}_{\mathbf{x}}^2  \mathbf{\gamma}_{\mathbf{x}}^2  \mathbf{\delta}_{\mathbf{x}}^2 \, .
\end{align}
and denote the induced rational map $T_d\, :\, \PP\, V(d) \to \PP\, V(8)$ by the same letter. 

\

We remark that it is hard to calculate the
values of $S_d$ (or $T_d$) on a general homogeneous form $f$ of degree $d$ without knowing the entire expression of $S_d$ (resp. $T_d$) as a polynomial
in the coefficients $A_{\mathbf{i}}$, which is awkward. One can, however, work directly with the symbolic expressions given in
(5) and (7) if one writes $f$ as a linear combination of $d$-th powers of linear forms:
\begin{gather}
f = \lambda_1 l_1^d + \dots + \lambda_N l_N^d, \quad \mathrm{some}\quad N\in\mathbb{N} \, .
\end{gather}
For linear forms $l_i$, $l_j$, $l_k$, $l_p\in \CC[x_1,\: x_2,\: x_3]_1$ we use the notation
\begin{gather}
I(l_i,\: l_j, \: l_k,\: l_p), \quad S_d(l_i,\: l_j, \: l_k,\: l_p) , \quad T_d(l_i,\: l_j, \: l_k,\: l_p)
\end{gather}
which is defined via formulas (4), (5), (7), but where for the vectors $\alpha$, $\beta$, $\gamma$, $\delta$ of symbolic variables we substitute the \emph{vectors of coordinates} w.r.t. $x_1$, $x_2$, $x_3$ of $l_i$, $l_j$, $l_k$, and $l_p$. One then has the following easy, but fundamental multi-linearity properties of $S_d$ and $T_d$ whose proof is a straight-forward computation and
therefore omitted.
\begin{lemma} 
\xlabel{lMultilinear} %HC
We have
\begin{align}
S_d (f ) = 24 \sum \lambda_i\lambda_j\lambda_k\lambda_p S_d ( l_i,\:  l_j,\:  l_k,\: 
l_p ) \, \\
T_d (f ) = 24 \sum \lambda_i\lambda_j\lambda_k\lambda_p T_d ( l_i,\:  l_j,\:  l_k,\: 
l_p ) 
\end{align}
The right-hand sums run over all $(i,\: j,\: k,\: p)$ with $1\le i< j<k<p\le N$.
\end{lemma}

\section{Special linear subspaces of the base loci}

The group $G=\mathrm{SL}_3\,\CC$ is a rank $2$ complex semisimple algebraic group, and choosing the standard torus $T$ of
diagonal matrices as maximal torus, and the group of upper-triangular matrices as Borel subgroup, one has the notions of
roots, positive and simple roots, and simple coroots $H_1$, $H_2$ available. Corresponding to $H_1$, $H_2$ one has
one-parameter subgroups $\lambda_{H_1},\:\lambda_{H_2}\, :\, \CC^{\ast} \to T$ given by
\begin{gather}
\lambda_{H_1} (t) = \left( \begin{array}{ccc} t & 0 & 0 \\ 0 & t^{-1} & 0 \\ 0 & 0 & 1 \end{array}\right) , \;
\lambda_{H_2} (t) = \left( \begin{array}{ccc} 1 & 0 & 0 \\ 0 & t & 0 \\ 0 & 0 & t^{-1} \end{array}\right) \, .
\end{gather}
For $d=3n+1$, we may view the covariant $S_d$ as an element in
\begin{gather}
\left( \mathrm{Sym}^4\, V(d)^{\vee} \otimes V(4) \right)^G ,
\end{gather}
a $G$-invariant polynomial of degree $4$ in the curve coefficients $A_{\mathbf{i}}$, $\mathbf{i} \in [d]^3$, $|\mathbf{i} |=d$,
with values in $V(4)$. As such it is a linear combination of monomials 
\begin{gather}
A_{\mathbf{i}}A_{\mathbf{j}}A_{\mathbf{k}}A_{\mathbf{l}} \mathbf{x}^{\mathbf{e}} ,
\end{gather}
where $\mathbf{i}$, $\mathbf{j}$, $\mathbf{k}$, $\mathbf{l}\in [d]^3$, $|\mathbf{i}|=\dots =|\mathbf{k}|=d$, $\mathbf{e}\in
[4]^3$, $|\mathbf{e}| = 4$.\\
Similarly, for $d=3n+2$, $T_d$ can be viewed as an element of 
\begin{gather}
\left( \mathrm{Sym}^4\, V(d)^{\vee} \otimes V(8) \right)^G ,
\end{gather}
i.e. a $G$-invariant polynomial of degree $4$ in the curve coefficients $A_{\mathbf{i}}$, $\mathbf{i} \in [d]^3$, $|\mathbf{i} |=d$,
with values in $V(8)$. It is a linear combination of monomials 
\begin{gather}
A_{\mathbf{i}}A_{\mathbf{j}}A_{\mathbf{k}}A_{\mathbf{l}} \mathbf{x}^{\mathbf{e}} ,
\end{gather}
where $\mathbf{i}$, $\mathbf{j}$, $\mathbf{k}$, $\mathbf{l}\in [d]^3$, $|\mathbf{i}|=\dots =|\mathbf{k}|=d$, $\mathbf{e}\in
[8]^3$, $|\mathbf{e}| = 8$.\\
The following proposition is an important ingredient in the proof of rationality. 

\begin{proposition} \xlabel{pBoe}
For $d=3n+1$, the projectivization of the linear space
\begin{gather}
%L_S = \PP\,\langle \mathbf{x}^{\mathbf{i}}= x_1^{i_1}x_2^{i_2}x_3^{i_3} \, |\, \mathbf{i}\in [d]^3, \: |\mathbf{i}|=d, \: i_1\ge 2n+3 \rangle_{\CC}
%\subset \PP\, V(d)
L_S = x_1^{2n+3} \cdot \CC[x_1,x_2,x_3]_{n-2} \subset V(d)
\end{gather}
is contained in the base scheme $B_S$ of the rational map $$S_d\, :\, \PP\, V(d) \dasharrow \PP\, V(4)$$ with a
full %HC
triple structure, i.e. 
\[
\mathcal{I}_{\PP (L_S )}^3 \supset \mathcal{I}_{B_S} \, .
\]
Similarly, for $d=3n+2$, the projectivization of the linear space
\begin{gather}
%L_T = \PP\,\langle \mathbf{x}^{\mathbf{i}}= x_1^{i_1}x_2^{i_2}x_3^{i_3} \, |\, \mathbf{i}\in [d]^3, \: |\mathbf{i}|=d, \: i_1\ge 2n+5 \rangle_{\CC}
%\subset \PP\, V(d)
L_S = x_1^{2n+5} \cdot \CC[x_1,x_2,x_3]_{n-3} \subset V(d)
\end{gather}
is contained in the base scheme $B_T$ of $$T_d\, :\, \PP\, V(d) \dasharrow \PP\, V(8)$$ with a 
full %HC
triple structure.
\end{proposition}
\begin{proof}
Regardless of whether $d=3n+1$ or $d=3n+2$, the conditions that the monomials in (14) or (16) are invariant under the actions of the one-parameter subgroups $\lambda_{H_1}$ resp.
$\lambda_{H_2}$ read
\begin{align}
i_1+j_1+k_1+l_1 - i_2-j_2-k_2-l_2 - e_1 +e_2 &=0 \quad \mathrm{resp.}\\
i_2+j_2+k_2+l_2 - i_3-j_3-k_3-l_3 - e_2 +e_3 &=0 \, . 
\end{align}
Now for $d=3n+1$ we get 
\begin{align}
4(3n+1) &= |\mathbf{i}|+|\mathbf{j}|+|\mathbf{k}|+|\mathbf{l}|\\
 &= 3 (i_2
+j_2+k_2+l_2) +(e_1-e_2) +(e_3-e_2)\nonumber\\
 &= 3 (i_2+j_2+k_2+l_2) + 4 - 3e_2 .\nonumber 
\end{align}
and for $d=3n+2$ one has 
\begin{align}
4(3n+2) &= |\mathbf{i}|+|\mathbf{j}|+|\mathbf{k}|+|\mathbf{l}|\\
 &= 3 (i_2
+j_2+k_2+l_2) +(e_1-e_2) +(e_3-e_2)\nonumber\\
 &= 3 (i_2+j_2+k_2+l_2) + 8 - 3e_2 .\nonumber 
\end{align}
In both cases then it follows that
\begin{align}
i_1+j_1+k_1+l_1 &= 4n + e_1 , \\
i_2+j_2+k_2+l_2 &= 4n + e_2 , \nonumber \\
i_3+j_3+k_3+l_3 &= 4n + e_3 .\nonumber
\end{align}
In particular, for $d=3n+1$, $i_1+j_1+k_1+l_1 \le 4n + 4$, which means that at most $1$ out of the $4$ indices $i_1$, $j_1$, $k_1$, $l_1$
 can be $\ge (4n+4)/2 +1 = 2n+3$. Since $\mathcal{I}_{L_S}$ is generated by those $A_{\mathbf{i}}$ with $i_1< 2n+3$, this proves the 
frist %HC
assertion.\\
For $d=3n+2$, $i_1+j_1+k_1+l_1 \le 4n + 8$, whence at most $1$ out of $i_1$, $j_1$, $k_1$, $l_1$ can be $\ge (4n+8)/2 +1=2n+ 5$, which proves the proposition.    
\end{proof}

\begin{remark} \xlabel{rLinearity}
By construction, $L_S$ (resp. $L_T$) have the following basic property: For $g\in V(d) \backslash L_S$ (resp. $g\in V(d) 
\backslash L_T$), the restriction $S_d\mid_{\PP (L_S+\CC g)}$ (resp. $T_d\mid_{\PP ( L_T+\CC g) }$) is linear.
\end{remark}

\section{Fiberwise surjectivity of the covariants} \xlabel{sFiberwise}

To begin with, we will show how some elements of $L_S$ (resp. $L_T$) can be written as sums of powers. For this %HC
let $K$ be a positive integer.
\begin{definition}
Let $\mathbf{b} = (b_1,\dots , b_K)\in\CC^K$ be given. Then 
we denote by %HC
\begin{gather}
p_i^{\mathbf{b}}(c) := \prod_{\stackrel{j\neq i}{1 \le j\le K}} \frac{c-b_j}{b_i-b_j}
\end{gather}
for $i=1,\dots , K$ 
% denote %HC
the interpolation polynomials of degree $K-1$ w.r.t. $\mathbf{b}$ in the one variable $c$. 
\end{definition}

\begin{lemma} \xlabel{lConstruction}
Let $\mathbf{b} = (b_1,\dots , b_K)\in\CC^K$, $b_i\neq b_j$ 
for $i\not=j$, and set %HC
% all $i$ and $j$, and put write %HC
$x=x_1$, $y=\lambda x_2 +\mu x_3$, $(\lambda,\: \mu )\neq (0,\: 0)$. Suppose $d > K$ and 
put $l_i:= b_i x +y$. Then for each $c\in\CC$ with $c\neq b_i$, $\forall i$,
\begin{gather}
f(c) = p_1^{\mathbf{b}}(c)l_1^d +\dots + p_K^{\mathbf{b}}(c) l_K^d - (cx + y)^d
\end{gather}
is nonzero and 
%is %HC
divisible by $x^K$.
\end{lemma}

\begin{proof}
The coefficient of the monomial $x^A y^B$ in $f(c)$  is equal to
\begin{gather*}
{ d \choose A } (p^{\mathbf{b}}_1(c) b_1^A + \dots + p^{\mathbf{b}}_{K} (c) b_{K}^A - c^A ). %HC
\end{gather*}
For %and for %HC 
$A\le K-1$ one has 
\[
c^A = p_1^{\mathbf{b}} (c) b_1^A + \dots + p^{\mathbf{b}}_{K} (c) b_{K}^A \, .
\]
for all $c$ %. %HC
by interpolation. %HC
\end{proof}

Choosing $K=2n+3$, we obtain elements $f(c)\in L_S$, and for $K=2n+5$ elements $f(c)\in L_T$. Now for $d=3n+1$ consider the diagram

%begin HC *\composite{{+}{\times}}
\xycenter{
	\PP(L_S+\CC g) \ar[d] \ar@{ *{ } }[r]|-*+{\subset}&\PP V(d) \ar@{-->}[r]^{S_d} \ar@{-->}[d]^{\pi_{L_S}}& \PP V(4)\\
	\bar{g} \ar@{ *{ } }[r]|-*+{\in}&\PP (V(d)/L_S)
	}

%\vspace{1cm}\\
%\setlength{\unitlength}{1cm}
%\begin{picture}(1,1.5)
%\put(3,-0.5){$\PP\, V(4)$}
%\put(2.5,1.2){$g\in \PP\, V(d)$}
%\put(4.5 , 1.2){$\dasharrow$}

%\put(3.5,0.4){\vector(0, -1){0.4}}
%\put(3.5, 0.7){\line(0, -1){0.2}}
%\put(3.5,1){\line(0, -1){0.2}}

%\put(5.5, 1.2){$\PP (V(d) / L_S)\ni \bar{g}=\pi_{L_S}(g)$}

%\put(4.6 , 1.7){${\pi}_{L_S }$}
%\put(3, 0.5){$S_d $}
%\end{picture}\vspace{1cm}\\

%end HC
or for $d=3n+2$ the diagram
%begin HC
\xycenter{
	\PP(L_T+\CC g) \ar[d] \ar@{ *{ } }[r]|-*+{\subset}&\PP V(d) \ar@{-->}[r]^{T_d} \ar@{-->}[d]^{\pi_{L_T}}& \PP V(8)\\
	\bar{g} \ar@{ *{ } }[r]|-*+{\in}&\PP (V(d)/L_T)
	}

%\vspace{1cm}\\
%\setlength{\unitlength}{1cm}
%\begin{picture}(1,1.5)
%\put(3,-0.5){$\PP\, V(8)$}
%\put(2.5,1.2){$g\in \PP\, V(d)$}
%\put(4.5 , 1.2){$\dasharrow$}

%\put(3.5,0.4){\vector(0, -1){0.4}}
%\put(3.5, 0.7){\line(0, -1){0.2}}
%\put(3.5,1){\line(0, -1){0.2}}

%\put(5.5, 1.2){$\PP (V(d) / L_T)\ni \bar{g}=\pi_{L_T}(g)$}

%\put(4.6 , 1.7){${\pi}_{L_T }$}
%\put(3, 0.5){$T_d $}
%\end{picture}\\
%\vspace{1cm}

%The basic genericity property we will have to establish in the following is this (case $d=3n+1$):\\ 
%{\bf Claim}: For a special (hence a general) $g$ the fibre $L_S + g$ of $\pi_L$ maps dominantly onto $\PP\, V(4)$ (similar assertion for $d=3n+2$).\\

The aim of this section is to prove

\begin{proposition} \xlabel{pSurjective}
Let $d=3n+1 \ge 37$. Then there exists a $g \in V(d)$ such that
\[
	S_d|_{\PP (L_S+\CC g)} \colon \PP (L_S+\CC g ) \dasharrow  \PP V(4)
\]
is surjective. For $d=3n+2 \ge 65$ there exists a $g \in V(d)$ such that
\[
	T_d|_{\PP (L_T+ \CC g)} \colon \PP ( L_T+\CC g ) \dasharrow \PP V(8)
\]
is surjective.
\end{proposition}

%We will now explain the method how this is done in the case $d=3n+1$. The case $d=3n+2$ is very similar, and we will deal with it afterwards.

We will prove the case $d=3n+1$ first. The case $d=3n+2$ is very similar, and we will deal with it afterwards.

We start by constructing points in the image $S_d$:

\begin{lemma} \xlabel{lCoefficients}
Consider $S_d (f(c) +g)$ as an element of $\CC[x_1,\: x_2,\: x_3, \: c]$ and write
\begin{gather}
S_d (f(c) +g) = Q_d c^d + \dots + Q_0
\end{gather}
with $Q_i\in\CC[x_1,\: x_2,\: x_3]_4$. Then $[Q_i ] \in S_d(\PP (L_S+\CC g))$ for all $i$.
\end{lemma}

\begin{proof}
The map
\begin{align}
\varphi \, :\, &\mathbb{A}^1 \to \PP\, V(4)\\
               &c \mapsto S_d (f(c) +g )\nonumber 
\end{align} 
gives a rational normal curve $X$ in $S_d (\PP (L_S +\CC g))$. Since by Remark 2.2, $S_d\mid_{\PP (L_S + \CC g)}$ is linear, the
linear span of
$X$
 is contained in $S_d (\PP (L_S +\CC g))$. Now $\langle X\rangle = \langle Q_0, \dots , Q_d \rangle$ which proves the claim.
\end{proof}
%\

%One may view $S_d (f(c) +g)$ as an element of $\CC[x_1,\: x_2,\: x_3, \: c]$ and write
%\begin{gather}
%S_d (f(c) +g) = Q_d c^d + \dots + Q_0, \quad Q_i\in\CC[x_1,\: x_2,\: x_3]_4\, .
%\end{gather}
%Now 
%\begin{align}
%\varphi \, :\, &\mathbb{A}^1 \to \PP\, V(4)\\
 %              &c \mapsto S_d (f(c) +g )\nonumber 
%\end{align} 
%gives a rational normal curve $X$ in $S_d (L_S +g)$. Since by remark 2.2, $S_d\mid_{L_S + g}$ is linear, $\langle X \rangle \subset S_d (L_S +g)$, so $\langle X\rangle = \langle Q_0, \dots , Q_d \rangle$, and the $Q_i$ are in $S_d (L_S + g)$.

Surprisingly, for $i$ large enough, the $Q_i$ do not depend on the vector $\mathbf{b} = (b_1,\dots , b_K)$ chosen to construct
$f(c)$:

\begin{proposition} \xlabel{pSurprise}
If 
\[
S_d (f(c) +g) = Q_d c^d + \dots + Q_0, \quad Q_i\in\CC[x_1,\: x_2,\: x_3]_4\, ,
\]
and 
\[
S_d ( -(cx + y)^d +g) = Q_d' c^d + \dots + Q_0', \quad Q_i' \in\CC[x_1,\: x_2,\: x_3]_4\, ,
\]
then $Q_i=Q_i'$ for $i\ge K$.
\end{proposition}

\begin{proof}
Write $g$ as a sum of $d$th powers of linear forms
\begin{gather}
g= m_1^d + \dots + m_{\mathrm{const}}^d\, ,
\end{gather}
where $\mathrm{const}$ is a positive integer that will be fixed 
(independently %(independent %HC
of $n$) in the later discussion. Then (using
%that %HC
$\mathcal{I}_{B_S}\subset \mathcal{I}_{\PP (L_S)}^3$ 
%) %HC 
and Lemma \ref{lMultilinear}) %HC
\begin{align}
S_d (f(c) + \epsilon g) &= S_d 
\Bigl( %HC
p_1^{\mathbf{b}}(c) l_1^d + \dots + p^{\mathbf{b}}_{K}(c) l_{K}^d - (cx + y)^d\\
                        &
                        \quad\quad%HC
                        + \epsilon m_1^d + \dots + \epsilon m_{\mathrm{const}}^d 
 			\Bigr) %HC
			 \nonumber \\
                       	&=24
                        	\Bigl( %\cdot \left\{  %HC
                        	\epsilon^3 \sum_{i< j < k < p}p^{\mathbf{b}}_i(c)I(l_i,\: m_j, \: m_k, \: m_p)^n  l_i m_jm_km_p 
			%\right. %HC
			\nonumber \\
	               	&
			\quad\quad %HC
			- \epsilon^3 \sum_{j< k< p} I(cx + y ,\: m_j,\: m_k,\: m_p)^n (cx  + y) m_j m_k m_p \nonumber \\
                        &\quad\quad %\left. %HC 
                        + \epsilon^4 \sum_{i < j < k < p}I(m_i,\: m_j, \: m_k, \: m_p)^n m_i m_j m_k m_p 
                        \Bigr) %\right\}  %HC
                        \nonumber 
\end{align}
For $\epsilon =1$ we find 
\begin{align}
S_d (f(c) +  g) &=  \sum_{i,\: j,\: k,\:  p}p^{\mathbf{b}}_i(c)I(l_i,\: m_j, \: m_k, \: m_p)^n  l_i m_jm_km_p\\
                        &
                        \quad %HC
                         + S_d( -(cx+y)^d + g) \, .\nonumber 
\end{align}
Since $\deg p_i^{\mathbf{b}} = K-1$, the assertion follows.
\end{proof}
Next we will investigate the dependence of $Q_t$ on $n$ for $t\ge K$.\\
We choose a fixed constant $\rmprec \in\mathbb{N}$ (the ``precision'') with $d-\rmprec > K$ and $\rmprec < n$ (later $\rmprec$ will be a prime number).  

\newcommand{\Fprec}{\mathbb{F}_{\rmprec}} %HC

\begin{lemma} \xlabel{lQuasipolynomials}
For $3n+1-\rmprec \le s \le 3n$, %$d-\rmprec \le s \le d$, %HC
the coefficient of $c^s$ in $I(cx + y, \: m_i,\: m_j, \: m_k )^n$ is of the form 
\begin{gather}
\varrho^n P(n)
\end{gather}
where $\varrho\in \CC$ (independent of $n$ and $s$) and $P(n)$ is a polynomial of degree 
$3n-s< \rmprec$. %$d-1-s< \rmprec$. %HC
$P(n)$ is, as a polynomial in $\CC[n]$ divisible by 
\[{ n \choose n-\lceil \frac{s}{3} \rceil }\, .\] %\[{ n \choose \lceil \frac{d-1-s}{3} \rceil }\, .\] %HC

%begin HC
If the coefficients of the $m_i$ are integers and $\rmprec$ is a prime number, then the reduction of
the coefficient of $c^s$ in $I(cx + y, \: m_i,\: m_j, \: m_k )^n$ modulo $\rmprec$ is still of the form 
\begin{gather}
\varrho^n P(n)
\end{gather}
whith $\varrho\in \Fprec$ and $P(n)\in \Fprec[n]$ satifying the same independence and divisibility
conditions as above.
%end HC

\end{lemma}

\begin{proof}
Calculating either over $\CC$ or over $\Fprec$, we have %We have %HC
\begin{gather}
I( cx+y,\: m_i,\: m_j,\: m_k)^n  \\
                                = (cx+y,\: m_i, \: m_k)^n (cx+y,\: m_i,\: m_j)^n (cx+y,\: m_j,\: m_k)^n (m_i,\: m_j,\: m_k)^n \nonumber \\
                                = (\xi_{ik} c +\eta_{ik} )^n (\xi_{ij} c +\eta_{ij})^n (\xi_{jk}c +\eta_{jk})^n (m_i,\:  m_j, \: m_k)^n\nonumber 
\end{gather}
where the $\xi$'s and $\eta$'s are constants 
%in $\CC$ %HC
(fixed once the $m$'s are fixed). If any of the $\xi$'s vanishes the polynomials $I( cx+y,\: m_i,\: m_j,\: m_k)^n$ is of 
degree $\le 2n$ in $c$. Since $s\ge 3n-\rmprec> 2n$, in this situation the coefficient of $c^s$ is $0$ and we are finished.
Assume therefore that the $\xi$'s are invertible. %HC

The above expression expands to 
%This expression is %HC
\begin{gather*}
 %= %HC
  \left( \sum_{p=1}^{n} {n\choose p} \xi_{ik}^p c^p \eta_{ik}^{n-p}  \right)\cdot  \left( \sum_{q=1}^{n} {n\choose q} \xi_{ij}^q c^q \eta_{ij}^{n-q}  \right)\\
\cdot \left( \sum_{r=1}^{n} {n\choose r} \xi_{jk}^r c^r \eta_{jk}^{n-r}  \right)\cdot (m_i,\: m_j,\: m_k)^n
\end{gather*}
and the coefficient of $c^s$ is %in this expression is %HC
\begin{gather*}
(m_i,\: m_j,\: m_k)^n \sum_{p+q+r = s} {n\choose p}{n\choose q}{n\choose r} \xi_{ik}^p\xi_{ij}^q\xi_{jk}^r \eta_{ik}^{n-p}\eta_{ij}^{n-q} \eta_{jk}^{n-r} \, .
\end{gather*}
Put $p' =n-p$, $q' = n-q$, $r' =n-r$ and rewrite this as %(recall that now $d=3n+1$) %HC
\begin{gather*}
\left( (m_i,\: m_j ,\: m_k) \xi_{ik} \xi_{ij} \xi_{jk} \right)^n \sum_{p' + q' +r' = 3n-s} %d-1-s} 
{n\choose p' }{n\choose q' }{n\choose r' } \xi_{ik}^{-p'} \xi_{ij}^{-q'} \xi_{jk}^{-r'} \eta_{ik}^{p' }\eta_{ij}^{q' } \eta_{jk}^{r' }
\end{gather*}
%(we can assume that $\xi_{ij}$, $\xi_{ik}$, $\xi_{jk}$ are invertible since otherwise $I( cx+y,\: m_i,\: m_j,\: m_k)^n$ is of 
degree $\le 2n$ in $c$,and since $s\ge d-\rmprec> K=2n+3$, the coefficient of $c^s$ is $0$ then). %HC
The first claim of the lemma 
over $\CC$ %HC
is obvious now. 
%begin HC
The reductions of the binomial coefficients modulo $\rmprec$ are polynomials in $n$ over $\Fprec$ if
$p',q',r' < \rmprec$. Our conditions on $s$ imply this, since
\[
	p',q',r' \le p'+q'+r' = 3n-s < \rmprec.
\]
%end HC
As for the stated divisibility property in $\CC[n]$ 
and $\Fprec[n]$%HC
, remark that in
\begin{gather*}
{n\choose p' }{n\choose q' }{n\choose r' }
\end{gather*} 
%with $p' + q' +r' = d-1-s$, at least one of $p'$, $q'$, $q'$ is $\ge\lceil \frac{d-1-s}{3} \rceil$. %HC
with $p' + q' +r' = 3n-s$, at least one of $p'$, $q'$, $q'$ is $\ge n-\lceil \frac{s}{3} \rceil$.  %HC
\end{proof}

\begin{proposition} \xlabel{pQuasipolynomials}
For $d-\rmprec+1 \le t \le d$, the coefficient of each monomial $\mathbf{x}^{\mathbf{i}}$ in $Q_t$ is of the form
\begin{gather}
\sum_{\nu = 1}^{ {\mathrm{const}\choose 3}} \varrho_{\nu }^n P_{\nu }(n)
\end{gather} 
where $\varrho_{\nu }\in\CC$ are constants (independent of $n$), and $P_{\nu }(n)$ are polynomials of degree $\le d-t < \rmprec$,
%$d-t \le \rmprec$,  which are divisible by 
%\[{ n \choose \lceil \frac{d-1-t}{3} \rceil }\, .\]
\[{ n \choose n-\lceil \frac{t}{3} \rceil }\, .\]

%begin HC
If $g$ can be written as sum of powers with integer coefficiens and $\rmprec$ is a prime number,
the same is true for the reduction of $Q_t$ mod $\rmprec$.
%end HC
\end{proposition}

\begin{proof}
$Q_t$ is the coefficient of $c^t$ in
\begin{gather}
(-24) \sum_{1\le i < \: j < \: k \le \mathrm{const}}  I(cx+y,\: m_i,\: m_j,\: m_k)^n (cx+y) m_i m_j m_k 
\end{gather}
(cf. (29)), so we may apply Lemma \ref{lQuasipolynomials} with $s=t$ and $s=t-1$.
\end{proof}

\begin{definition} \xlabel{dR}
For $d-\rmprec+1 \le t \le d$, we put
\begin{gather}
%R_t := \frac{ Q_t}{{ n \choose \lceil \frac{d-1-t}{3} \rceil }}\, .
R_t := \frac{ Q_t}{{ n \choose n-\lceil \frac{t}{3} \rceil }}\, .
\end{gather}
\end{definition}

%begin HC
\begin{proof}[Proof of Proposition \ref{pSurjective} (for $d=3n+1$)]
Let $\rmconst=9$ and consider
\[
	g = m_1^d + \dots + m_\rmconst^d
\]
with 
{\tiny
\begin{align*}
	m_1 &= x_1+3 x_2+9 x_3
	&m_4 &= x_1+6 x_2-10 x_3
	&m_7 &= -3 x_2+2 x_3\\
	m_2 &= -10 x_1+x_2+4 x_3
	&m_5 &= 4 x_1-8 x_2-10 x_3	
	&m_8 &= 8 x_1-4 x_2-4 x_3\\
	m_3 &= 8 x_1+4 x_2+6 x_3
	&m_6 &= -3 x_1+7 x_2-4 x_3
	&m_9 &= -10 x_1+4 x_2+6 x_3.
\end{align*}}

For $\rmprec=11$ we perform our construction with 
$x=x_1$ and two different values for $y$, namely $y_1=x_2$ and $y_2=x_3$. 
We obtain $22$ quartics
$R_d^{y_1}, \dots , R^{y_1}_{d-10}, 
R_d^{y_2}, \dots , R^{y_2}_{d-10}$.
By Lemma \ref{lCoefficients} and Proposition \ref{pSurprise} these quartics are in the image of
$S_d|_{L_S+g}$ if 
$$d-10 \ge K  \iff 3n+1-10 \ge 2n+3 \iff n \ge 12.$$
The coefficients of the $R^{y_i}_j$ form a $15 \times 22$ matrix $M(n)$ with entries of the form
$\sum_{\nu = 1}^{84} \varrho_{\nu }^n P_{\nu }(n)$ by Proposition \ref{pQuasipolynomials}. 
Modulo $11$ this matrix becomes
periodic in $n$ with period $11\cdot10 = 110$. With a computer algebra program it is
straightforward to check that all these matrices have full rank $15$. A {\ttfamily Macaulay2} script doing this can be found at
\cite{rationalityScript}. This proves the claim for $d=3n+1$.
\end{proof}

%Hence, for $d= 3n+1$, we can establish the Claim on page 6 via the following algorithm using a computer algebra system,e.g.
%{\ttfamily Macaulay 2}:

%\

%\flushleft {\bf Algorithm} (case $d=3n+1$): We choose $\mathrm{const}=9$ and $x=x_1$ and for $y$ the three different values $y_1=x_2+x_3$, $y_2=x_2$, $y_3=x_3$. We choose $\rmprec=17$. Finally we select $m_1, \dots , m_9\in\mathbb{Z}[x_1,\: x_2,\: x_3]_1$ in such a way that in the polynomials
%\begin{gather*}
%R_d^{y_1}, \dots , R^{y_1}_{d-\rmprec+1}, \\
%R_d^{y_2}, \dots , R^{y_2}_{d-\rmprec+1}, \\
%R_d^{y_3}, \dots , R^{y_3}_{d-\rmprec+1}, \\
%\end{gather*}
%(we get one set of $\rmprec$ polynomials via definition 3.6 for each choice of $y$) all prime factors of occurring denominators are $< \rmprec$. Remark that these polynomials are in $\mathbb{Q}[x_1, \: x_2,\: x_3]_4$ here. We give the explicit choices for the $m$'s in the Appendix.\\
%We reduce the $R$'s mod $\rmprec$, and notice that, thanks to proposition 3.5, the $(3\cdot \rmprec)\times 15$ matrix consisting of the coefficients of all the monomials $\mathbf{x}^{\mathbf{i}}$ in the $(3\cdot \rmprec)$ quartics given by the $R$'s, is \emph{periodic} as a matrix of functions of $n$ \emph{with period} $\rmprec\cdot (\rmprec -1)$.\\
%Hence we check on a computer that all of these $\rmprec\cdot (\rmprec -1)$ matrices over $\mathbb{F}_{\rmprec}$ are of full rank.

%\

Let us turn to the case $d=3n+2$. The whole procedure is similar in this case. If we take $K = 2n+5$  
Lemma \ref{lConstruction}, Proposition \ref{pSurprise}, Lemma \ref{lQuasipolynomials} and Proposition \ref{pQuasipolynomials} remain true as stated and Definition \ref{dR} still makes sense.

\begin{proof}[Proof of Proposition \ref{pSurjective} (for $d=3n+2$)]
Let $\rmconst=9$ and consider
\[
	g = m_1^d + \dots + m_\rmconst^d
\]
with $m_i$ as above.

For $\rmprec=19$ we perform our construction with 
$x=x_1$ and three different values for $y$, namely $y_1=x_2$, $y_2=x_3$ and $y_3=x_2+x_3$.
We obtain $57$ octics
$R_d^{y_1}, \dots , R^{y_1}_{d-18}, 
R_d^{y_2}, \dots , R^{y_2}_{d-18},
R_d^{y_3}, \dots , R^{y_3}_{d-18}$.
By Lemma \ref{lCoefficients} and Proposition \ref{pSurprise} these octics are in the image of
$S_d|_{L_T+g}$ if 
$$d-18 \ge K  \iff 3n+2-18 \ge 2n+5 \iff n \ge 21.$$
The coefficients of the $R^{y_i}_j$ from a $45 \times 57$ matrix $M(n)$ with entries of the form
$\sum_{\nu = 1}^{84} \varrho_{\nu }^n P_{\nu }(n)$ by Proposition \ref{pQuasipolynomials}. 
Modulo $19$ this matrix becomes
periodic in $n$ with period $19\cdot18 = 342$. With a computer algebra program it is
straightforward to check that all these matrices have full rank $45$. A {\ttfamily Macaulay2} script doing this can be found at
\cite{rationalityScript}. This proves the claim for $d=3n+2$.
\end{proof}

%end HC

\newcommand{\verystable}{\bigl(\PP\, V(4)\bigr)_{\mathrm{vs}}}

\section{Sections of principal bundles and proof of rationality}

We will now show how to conclude the proof in the case $d=3n+1$. We make some comments on the case $d=3n+2$ when they are in order, but otherwise leave the obvious modifications to the reader.
Let $\verystable \subset\PP\, V(4)$ be the open subset of very stable points with respect to the
action of
$\bar{G}$ and the
$\bar{G}$-linearized line bundle $\mathcal{O} (3)$ (very stable means stable with trivial stabilizer). Now the essential point
is

\begin{proposition}
The quotient morphism
\[
\verystable \to \verystable / \bar{G}
\]
is a principal $\bar{G}$-bundle in the Zariski topology. 
\end{proposition}
\begin{proof}
See \cite{Shep}, Prop. 2. This holds also true with $V(4)$ replaced with $V(8)$.
\end{proof}
It follows that this $\bar{G}$-bundle has a section defined generically which we will denote by $\sigma_4$.

\begin{proof}[Proof of Theorem \ref{tMain}]
Consider the incidence variety
\[
	X = \{ (g,\bar{g},f) \,|\, \pi_{L_S}(g)=\bar{g}, S_d(g)=f\} 
	\subset \PP V(d) \times \PP\bigl(V(d)/L_S\bigr) \times \PP V(4)	
\]
and the diagram
\newcommand{\rmpr}{\mathrm{pr}}
\xycenter{
	X \ar[d]_{\rmpr_{23}} \ar@{<-->}[r]^{1:1}_{\rmpr_{1}}
	& \PP V(d) \ar@{-->}[r]
	& \PP V(d)/\bar{G} \ar@{-->}[dd]^{\bar{S}_d}
	\\
	\PP\bigl( V(d)/L_S \bigr) \times \PP V(4) \ar[d]
	\\
	\PP V(4) \ar@{-->}[rr]
	& & \PP V(4)/\bar{G} \ar@/^20pt/@{-->}[ll]^{\sigma_4}.
	}
By Proposition \ref{pSurjective} the projection $\rmpr_{23}$ is dominant. It follows then from Remark \ref{rLinearity} that $X$ is birational to a vector bundle over $\PP\bigl( V(d)/L_S \bigr) \times \PP V(4)$ and hence also over $\PP V(4)$. After replacing $\sigma_4$ by a translate, we can assume that $\sigma_4$ meets an open set $U \subset \PP V(4)$ over which $X$ is trivial.  Since $\bar{G}$ acts generically freely on $\PP V(4)$, we can pull back the above vector bundle structure via $\sigma_4$ and obtain that 
$\PP V(d)/\bar{G}$ is birational to $\PP V(4)/\bar{G} \times \PP^N$ with $N = \dim V(d) - \dim V(4)$. If $d \ge 37$ as in Proposition
\ref{pSurjective}, then certainly $N \ge 8$ and  since $\left( \PP V(4) /\bar{G} \right) \times \PP^8$ is rational, $\PP V(d)/\bar{G}$ is rational.  
In the case $d=3n+2$ the same argument works since the space of octics is also stably rational of level $8$. This proves Theorem \ref{tMain}.
\end{proof}

\end{document}